\documentclass[11pt]{amsart} \usepackage{latexsym}
\usepackage{amsfonts}

\usepackage{amsmath,amsthm,amssymb}

\newtheorem{theor}{Theorem} 
 
\newtheorem{thm}{Theorem}[section] \newtheorem{lem}[thm]{Lemma}
\newtheorem{obs}[thm]{Observation} 
\newtheorem{conj}[thm]{Conjecture}
\newtheorem{prop}[thm]{Proposition} \theoremstyle{definition}
\newtheorem{defn}[thm]{Definition}

\newtheorem{conv}[thm]{Convention} \newtheorem{rem}[thm]{Remark}

\begin{document}

\title[Delzant's $T$-invariant and Kolmogorov Complexity] {Delzant's
  $T$-invariant, Kolmogorov complexity and one-relator groups }

\author[I.~Kapovich]{Ilya Kapovich}

\address{\tt Department of Mathematics, University of Illinois at
  Urbana-Champaign, 1409 West Green Street, Urbana, IL 61801, USA
  \newline http://www.math.uiuc.edu/\~{}kapovich/} \email{\tt
  kapovich@math.uiuc.edu}

\author[P.~Schupp]{Paul Schupp}

\address{\tt Department of Mathematics, University of Illinois at
  Urbana-Champaign, 1409 West Green Street, Urbana, IL 61801, USA
  \newline http://www.math.uiuc.edu/People/schupp.html}
\email{schupp@math.}
\begin{abstract}
  We prove that for ``random'' one-relator groups the Delzant
  $T$-invariant (which measures the smallest size of a finite
  presentation of  a group) is comparable in magnitude with the length
  of the defining relator. The proof relies on our previous results
  regarding isomorphism rigidity of generic one-relator groups and on
  the methods of the theory of Kolmogorov-Chaitin complexity.
\end{abstract}

\thanks{The authors
  were supported by the NSF grant DMS\#0404991 and the NSA grant
  DMA\#H98230-04-1-0115.}


\maketitle


\section{Introduction}

Delzant~\cite{De} introduced an extremely interesting but still rather
enigmatic invariant for finitely presentable groups.  For any finite
presentation
\newline $\Pi=\langle X | R\rangle$ define \emph{the length}
$\ell(\Pi)$ as
\[
\ell(\Pi):=\sum_{r\in R} \max\{ |r|-2, 0\}.
\]
If $G$ is a finitely presentable group, the \emph{T-invariant $T(G)$}
of $G$, which we  also call the \emph{presentation rank of $G$}, is
defined~\cite{De} as
\[
T(G):=\min\{\ell(\Pi) | \Pi \text{ is a finite presentation of the
  group} ~G\}.
\]

The $T$-invariant plays a central role in Delzant and Potyagailo's
proof of the strong accessibility (or "hierarchical decomposition")
theorem for finitely presented groups~\cite{DP}.  This theorem is the
strongest and most difficult of numerous accessibility
results~\cite{Du85,Du93,BF91,BF91a,Se1,De1,W,KW}.  One can also define
a closely related notion, the \emph{non-reduced $T$-invariant
  $T_1(G)$}, as the minimum of sums of lengths of the defining
relators, taken over all finite presentations of $G$. As we observe in
Lemma~\ref{lem:l_1} below, if $G$ is a finitely presentable group
without elements of order two then
\[
T(G)\le T_1(G)\le 3T(G).
\]
The non-reduced $T$-invariant has been studied in the context of
3-manifolds, where it turns out to be related to the notion of Matveev
complexity. We refer the reader to a recent paper of Pervova and
Petronio~\cite{PP} for a discussion on this subject.

If $G$ is a finitely generated group then the ordinary  \emph{rank}, $rk(G)$, of
$G$ is the smallest cardinality of a finite generating set for $G$.
The first (and already quite nontrivial) accessibility result
is Grushko's theorem~\cite{Gru} which asserts that for finitely
generated groups $G_1$ and $G_2$ we have $rk(G_1\ast
G_2)=rk(G_1)+rk(G_2)$. In \cite{De} Delzant proved a similar theorem
for the presentation rank, namely that
\[
T(G_1\ast G_2)=T(G_1)+T(G_2)
\]
if $G_1,G_2$ are finitely presentable groups.

The hierarchical decomposition theorem proved in \cite{DP} implies,
for example, that an iterated process of JSJ-decomposition (in any
sense of the word) \cite{Se2,RiSe,DuSa,FuPa,Bow98} applied to a
finitely presented group, then to the factors of its
JSJ-decomposition, and so on, always terminates. The $T$-invariant is
also crucial in Delzant's generalization~\cite{De1} of Sela's
acylindrical accessibility result~\cite{Se1} for finitely presented
groups.

If $\Pi$ is a finite presentation, let $G(\Pi)$ be the group
defined by $\Pi$. We can regard $T$ as a function defined over
finite presentations by setting
\newline $T(\Pi) = T(G(\Pi))$. If $G$ is
given by a particular finite presentation $\Pi$ then $\ell(\Pi)$
gives an obvious upper bound for $T(G(\Pi)))$. However,
it is very unclear  in general how to estimate $T(G)$ from below. For example,
if $\Pi=\langle X| R\rangle$ and $\alpha\in Aut(F(X))$ then the
presentations $\Pi$ and $\Pi'=\langle X | \alpha(R)\rangle$ define
isomorphic groups but it is easy to produce examples where
$\ell(\Pi')$ is arbitrarily smaller than $\ell(\Pi)$.

We prove however that for "most" one-relator presentations this does not
happen and that the value of Delzant's $T$-invariant is comparable in
magnitude with the length of the defining relator. If $r\in
F(a_1,\dots, a_k)$,
\newline let $G_r:=\langle a_1,\dots, a_k | r\rangle$ be
the one-relator group whose defining relator is $r$.  Our main result
is:

\begin{theor}\label{A}
   Fix an  integer  $k>1$ and let $F=F(a_1,\dots, a_k)$.  For any
  number $0<\epsilon<1$ there is an integer $n_1>0$ and a constant
  $M=M(k,\epsilon)>0$ with the following property.

  Let $J$ be the set of all nontrivial cyclically reduced words $r$
  such that
\[
T(G_r)\log_2 T(G_r) \ge M |r|.
\]
Then for any $n\ge n_1$
\[
\frac{\#\{r\in J : |r|=n \}}{\#\{r\in F : r \text{ is cyclically
    reduced and } |r|=n \}}\ge 1-\epsilon.
\]
\end{theor}

Thus for any fixed $0<\epsilon, \delta<1$ we asymptotically have
$T(G_r)\ge c |r|^{1-\delta}$, where $c$ is a constant, for at least
the fraction $(1-\epsilon)$ of all cyclically reduced words $r$ of a
given length. This says that the description of a one-relator group by
a generic relator $r$ is ``essentially incompressible''.  In view of the above
remarks about the connection between $T(G)$ and $T_1(G)$, the same
conclusion as in Theorem~\ref{A} also holds for $T_1(G_r)$.

This is a good place to observe that the function $T$ is not
computable.

\begin{obs}\label{unsolvable}
  The function $T$, as a function over finite presentations, is not a
  computable function.
\end{obs}

\begin{proof}
  We say that a finitely generated group $G$ is \emph{essentially
    free} if $G$ is the free product of a finitely generated free
  group and finitely many cyclic groups of order two. The only
  defining relators in the ``standard presentation'' $\Pi_0$ of such a
  group are the squares of those generators which have order two and
  so $\ell(\Pi_0)=0$.

  It is easy to use Tietze transformations to show that any group $G$
  having a finite presentation in which all relators have length at
  most two is essentially free. Hence, by the definition of $T(G)$, a
  finitely presentable group $G$ has $T(G) = 0$ if and only if $G$ is
  essentially free.

  Recall that a property $\mathcal{P}$ of finitely presented groups is
  a \emph{Markov property} if $\mathcal{P}$ is independent of
  presentation, there are finitely presented groups with $\mathcal{P}$
  and there is a finitely presented group $G_*$ such that $G_*$ cannot
  be embedded in any finitely presented group with $\mathcal{P}$.
  Being essentially free is clearly a Markov property.  We can take
  $G_*$ to be the cyclic group of order three.  The classic
  Adian-Rabin Theorem \cite{LS} says that if $\mathcal{P}$ is any
  Markov property then there is no algorithm over all finite
  presentations which, when given a finite presentation $\Pi$, decides
  whether or not the group $G(\Pi)$ has $\mathcal{P}$.

  If the function $T$ were computable then, for  any finite
  presentation $\Pi$, we could decide the  essential freeness of   $G(\Pi)$
  by computing $T(\Pi)$.  Hence $T$ cannot be computable.
\end{proof}

The proof of Theorem A involves several different probabilistic tools.
The  idea introduced in this paper is the use of Kolmogorov
complexity, a concept that plays an important role in coding theory,
algorithmic probability and complexity theory. This notion is also
sometimes known as ``Kolmogorov-Chaitin complexity'' because of the
contributions  of Chaitin to the subject.  Roughly speaking, the
Kolmogorov complexity of a word is the size of the smallest computer
program (in a fixed programming language) that can compute the given word.
Surprisingly, the only previous use of Kolmogorov complexity in group theory
known to us is a 1985 paper of Grigorchuk~\cite{Gri}, giving an  interesting
application of Kolmogorov complexity to algorithmic
problems in  group theory.

  Our results here also depend  on  \cite{KS} and \cite{KSS}
where we obtained a number of results
regarding a very strong Mostow-type "isomorphism rigidity" for generic
one-relator groups. These  results use  a combination of the
Arzhantseva-Ol'shaskii minimization technique and their ingenious
"non-readability" small cancellation condition \cite{KS} and
Large Deviation Theory  \cite{KSS} to study  the behavior of  random words
under an arbitrary automorphism of the ambient free group.
The isomorphism rigidity theorems proved in \cite{KSS}
allow us, given any finite presentation $\Pi=\langle X |R \rangle$
defining a group isomorphic to $G_r=\langle a_1, \dots, a_k|
r\rangle$ (where $k>1$ is fixed) for a generic relator $r$ plus a
small initial segment $u$ of $r$, to algorithmically recover the
word $r$. This implies that $r$ is uniquely algorithmically
determined by an amount $O(\ell(\Pi)\log \ell(\Pi))$ of
information. (The logarithmic term comes from the fact that the
subscripts in the enumeration of letters in $X$ also need to be
encoded.) From here one can deduce that the Kolmogorov
complexity of $r$ is $\le O(\ell(\Pi)\log \ell(\Pi))$. On the
other hand, using the methods of algorithmic probability, in
particular the notion of \emph{prefix complexity}, we can deduce
that a cyclically reduced word $r$ of a given length has
Kolmogorov complexity $\ge c |r|$ asymptotically with probability
$\ge 1-\epsilon$. These inequalities  taken together yield the
conclusion of Theorem~\ref{A}.

We believe that the general analogue of Theorem~\ref{A} is true. This
would say that if we fix a number $k\ge 2$ of generators and any
number $m\ge 1$ of defining relators, then a generic $k$-generator
$m$-relator presentation should essentially be the shortest
description of the group defined.  We have seen that the proof of
Theorem~\ref{A} relies on two components: the Kolmogorov complexity
arguments used in this paper and the isomorphism rigidity results for
random one-relator groups established in ~\cite{KS,KSS}.  Now most of
the arguments and statements of~\cite{KS,KSS} needed to prove
isomorphism rigidity actually go through for generic groups with an
arbitrary fixed number of relators and we believe that ``generic groups
are rigid'' in general.  However, to actually infer rigidity, at the
end of the proof we use a crucial fact about one-relator groups.
Namely, we need the classical theorem of Magnus (see, for example,
\cite{LS}) which says that if two elements $r$ and $s$ have the same
normal closures in a free group $F$ then $r$ is conjugate in $F$ to
$s$ or $s^{-1}$.  This statement does not hold in general for tuples
consisting of more than one element of $F$. However, we believe that
the desired analogue does hold generically.

If $\tau = (u_1, \dots ,u_m)$ is an $m$-tuple of elements of the free
group $F_k$, the \emph{symmetrized set} $R(\tau)$ \emph{generated by }
$\tau$ consists of all the cyclic permutations of cyclically reduced
forms of $u_i^{\pm 1}$.

\begin{conj}[Stability  Conjecture]
  Fix $k\ge 2$ and $m \ge 1$ and let $F=F(a_1,\dots, a_k)$. Then there
  exists an algorithmically recognizable generic class $\mathcal{C}$
  of $m$-tuples of elements of $F$ with the following property. If
  $\sigma,\tau \in \mathcal{C}$ and $\alpha\in Aut(F)$ are such that
  $R(\sigma)$ and $R(\alpha(\tau))$ have the same normal closure in
  $F$ then $R(\sigma) = R(\alpha(\tau))$.
\end{conj}
Magnus' theorem implies that the Stability Conjecture holds for $m=1$
with $\mathcal C=F_k$. If one could establish the Stability
Conjecture, then both the isomorphism rigidity results of \cite{KSS}
and the results of this paper would then follow for finitely presented
groups with any fixed numbers of generators and relators exactly as in
the one-relator case.

In~\cite{KSS} we showed that for a fixed $k\ge 2$ the number $I_n$ of
\emph{isomorphism types} of $k$-generator one-relator groups with
cyclically reduced defining relators of length $n$ satisfies
\[
\frac{c_1(2k-1)^n}{n}\le I_n\le \frac{c_2(2k-1)^n}{n},
\]
where $c_1=c_1(k)>0, c_2=c_2(k)>0$ are some constants independent of
$n$. Using auxiliary results from the proof of Theorem~\ref{A} we
obtain an improvement of this estimate in the present paper and
compute the precise asymptotics of $I_n$:

\begin{theor}\label{B}
  Let $k\ge 2$ be a fixed integer. Then the number $I_n$ of
  isomorphism types of $k$-generator one-relator groups with
  cyclically reduced defining relators of length $n$ satisfies:

\[
I_n\sim \frac{(2k-1)^n}{nk!2^{k+1}}.
\]
\end{theor}
Here $f(n)\sim g(n)$ means that $\displaystyle\lim_{n\to\infty}
f(n)/g(n)=1$.

The authors are grateful to Carl Jockusch and Paul Vitanyi for helpful
discussions regarding Kolmogorov complexity. They also thank Warren
Dicks for suggesting the problem of  computing  the precise asymptotics of $I_n$.

\section{Kolmogorov Complexity}

The $T$-invariant is a measure of ``smallest descriptive complexity''
in the framework of finite presentations of groups while  Kolmogorov
complexity is a  general theory of ``minimal descriptive
complexity''.  We provide here only a brief discussion of the relevant
facts regarding Kolmogorov complexity and refer the reader to the
survey of Fortnow~\cite{Fort} for an overview and to the excellent and
comprehensive book of Li and Vitanyi~\cite{LV} for detailed background
information.

Intuitively speaking, the Kolmogorov complexity $C(x)$ of a finite
binary string $x$ is the size of the smallest computer program $M$
that can compute $x$. In order for this notion to make sense one
needs to first fix a ``programming language'' but it turns out that
all reasonable choices yield measures which are equivalent up to an
additive constant.

Note that  $C(x)$, as a measure of \emph{descriptive} complexity of $x$,
totally disregards how long the particular program $M$ will have to
run in order to compute $x$. Some strings clearly admit much shorter
descriptions then their length. For example, if $x_0$ is the binary
representation of the number $2^{2^{2^{10}}}$ then the length of $x_0$
is huge, namely $1+2^{2^{10}}$. Yet we were just able to give a very
short unambiguous description of $x_0$. Thus $x_0$ has small
Kolmogorov complexity and $C(x_0)<< |x_0|$.  On the other hand it is
intuitively clear that for a ``random'' string $x$ of large length,
the shortest description of $x$ is essentially $x$ itself. In this
case $C(x)\approx |x|$. This phenomenon is called
``incompressibility'' and plays an important role in complexity theory
for establishing  lower complexity bounds of various
algorithms.

Recall that any  Turing machine $M$ on the set of finite binary strings
$\{0,1\}^{\ast}$ computes a partial recursive function
$\{0,1\}^{\ast}\to \{0,1\}^{\ast}$ and, moreover, every partial
recursive function $\{0,1\}^{\ast}\to \{0,1\}^{\ast}$ arises in this
fashion.

Once one has fixed the formalism of Turing machines, one can identify
a Turing machine with its sequence of instructions and think of Turing
machines as programs.  A Turing machine $M$ can then itself be coded
as a binary string according to some fixed effective method
and we write $\langle M \rangle$ for the code of the machine $M$. The pair
consisting of a Turing machine $M$ and an input $w$ can then be given
the  code $\langle M\rangle w$.  A basic feature of the theory of
computability is the existence of a \emph{universal} Turing machine
$U$, which, if its input is a code $\langle M\rangle w$, simulates $M$
on input $w$. To be more precise, a Turing machine $U$ is
\emph{universal} if for any Turing machine $M$ there is a binary
string $\langle M\rangle$ such that for any string $w\in
\{0,1\}^{\ast}$ the machine $U$ produces the same result on input
$\langle M\rangle w$ as $M$ does on $w$.

\begin{defn}\label{defn:kolm}
  Fix a universal Turing machine $U$ with the alphabet
  \newline $\Sigma:=\{0,1\}$. Then  $U$ computes a universal partial recursive
  function $\phi$ from $\Sigma^{\ast}$ to $\Sigma^{\ast}$. That is,
  for any partial recursive function $\psi$ there is a string
 $z\in \Sigma^{\ast}$ such that for all  $x\in \Sigma^{\ast}$,
  $\phi (zx) = \psi(x)$.

  For a finite binary string $x\in \Sigma^{\ast}$ we define the
  \emph{Kolmogorov complexity} $C(x)$ as

\[
C(x):=\min \{ |p| : p\in \Sigma^{\ast}, \phi(p)=x\}.
\]

\end{defn}

Kolmogorov complexity is traditionally defined for finite binary
strings. However, if $s>1$ is a fixed integer, then all of the
standard definitions and theorems go through essentially unchanged if
one considers finite strings $x$ in a fixed $s$-letter alphabet $A$.
This can be done in either of two essentially equivalent ways. First,
one can modify Definition~\ref{defn:kolm} by choosing $U$ to be a
universal Turing machine with the alphabet $A_s$ computing a universal
partial recursive function from $A_s^{\ast}$ to $A_s^{\ast}$.
Alternatively, one can fix a recursive bijection $h: A_s^{\ast}\to
\Sigma^{\ast}$ and define $C_s(x)$, where $x\in A_s^{\ast}$ to be
$C(h(x))$. We choose the latter option since most theorems in
\cite{LV} are stated for binary strings and
we want to be able to cite the results of \cite{LV} verbatim.

\begin{defn}\label{defn:kolm1}
  Let $s>1$ be an integer and let $A_s$ be an alphabet with $s$
  letters.  Fix a recursive bijection $h: A_s^{\ast}\to
  \{0,1\}^{\ast}$.

  For any string $x\in A_s^{\ast}$ define its \emph{Kolmogorov
    complexity} $C_s(x)$ as

\[
C_s(x):=C(h(x)).
\]
\end{defn}

Kolmogorov complexity lacks some mathematical properties which are
essential for certain arguments.  Fortunately, this difficulty can be
overcome by using the closely related notion of \emph{prefix
  complexity}.  For a detailed discussion of this notion we refer the
reader to Chapters 2 and 3 of \cite{LV}. In the present paper we need
only cite a few basic facts regarding prefix complexity from
\cite{LV}.  A partial recursive function $\phi$ on $\Sigma^{\ast}$ is
called a \emph{prefix function} if whenever $\phi(x)$ is defined and
$x$ is a proper initial segment of $y$, then $\phi(y)$ is undefined.
There is a corresponding notion of a \emph{prefix machine}. Informally
speaking, a prefix machine does not require an ``end-of-tape'' symbol
for the input word and decides whether or not to halt only based on
its current state and before scanning the next letter of the input.
The machine starts working on an infinite input word and, after
performing a computational step on the working and output tapes, the
machine either moves one letter to the right on the input tape or
halts and terminates its work.

Just as with ordinary Turing machines, there exist universal prefix
machines computing universal prefix partial recursive functions (see
Theorem 3.1.1 in \cite{LV}).

\begin{defn}\label{defn:prefix}
  Fix a universal prefix Turing machine $U'$ with the alphabet
  $\Sigma=\{0,1\}$. Then  $U'$ computes a universal prefix partial
  recursive function $\psi$ from $\Sigma^{\ast}$ to $\Sigma^{\ast}$.

  For a finite binary string $x\in \Sigma^{\ast}$ we define the
  \emph{prefix complexity} $K(x)$ as

\[
K(x):=\min \{ |p| : p\in \Sigma^{\ast}, \psi(p)=x\}.
\]
\end{defn}

Similarly to the case of Kolmogorov complexity, prefix complexity can
be defined not only for binary but also for $s$-ary strings.

\begin{defn}\label{defn:prefix1}

  Let $s>1$ be an integer and let $A_s$ be an alphabet with $s$
  letters.  Fix the same recursive bijection $h: A_s^{\ast}\to
  \{0,1\}^{\ast}$ as in Definition~\ref{defn:kolm1}.

  For any string $x\in A_s^{\ast}$ define its \emph{prefix complexity}
  $K_s(x)$ as
\[
K_s(x):=K(h(x)).
\]

\end{defn}

For our purposes, the crucial  way in which prefix complexity
is  better than Kolmogorov complexity is that that $\sum_{x\in
  \{0,1\}^{\ast}} 2^{-K(x)}\le 1$ while $\sum_{x\in \{0,1\}^{\ast}}
2^{-C(x)}$ diverges.

We list here some relevant properties of Kolmogorov and prefix
complexity.

\begin{prop} \label{prop:use}

  Let $s>1$ be a fixed integer and let $A_s$ be an $s$-letter
  alphabet.  Then:

\begin{enumerate}
\item We have \[ \sum_{x\in\{0,1\}^{\ast}} 2^{-K(x)}\le 1.
\]
\item Up to additive constants for any $x\in \{0,1\}^{\ast}$ we have
\[
C(x)\le K(x)\le C(x)+\log_2 C(x).
\]
\item We have
\[
\sum_{x\in A_s^{\ast}} 2^{-K_s(x)}\le 1.
\]
\item Up to additive constants for any $x\in A_s^{\ast}$ we have
\[
C_s(x)\le K_s(x)\le C_s(x)+\log_2 C_s(x).
\]
\end{enumerate}

\end{prop}

\begin{proof}

  Part (1), as observed by Levin~\cite{Le}, is a direct corollary of
  Kraft's Inequality, which is ubiquitous in information theory (see
  also 4.2.2(b) in \cite{LV}). Part (2) is statement 3.1.3 in
  \cite{LV}. Clearly, (1) implies (3) and, also, (2) implies (4).
  Since part (1) is quite important for our purposes, we provide a
  proof here.

  A subset $S\subseteq \{0,1\}^{\ast}$ is \emph{prefix-free} if
  whenever $p,q\in S, p\ne q$ then $p$ is not an initial segment of
  $q$. Recall that by definition $K(x)$ is the shortest length of a
  prefix program $p$ with $\psi(p) = x$.  Thus the set $S$ of such all such
  $p$ corresponding to $ x\in \{ 0,1 \}^{\ast}$ is prefix-free.  If $p$ is a
  binary string, then $2^{-|p|}$ is the Lebesque measure of the subset
  $S_p$ of the unit interval $I = [0,1]$ consisting of those numbers
  whose binary expansion begins with $p$.  Since $S$ is prefix-free,
  subsets $S_p$ and $S_q$ are disjoint for $p \ne q$.  The inequality
  thus follows from the countable additivity of Lebesgue measure.

\end{proof}

We also recall the classical Markov inequality from probability theory
which can be found in most probability textbooks (see, for example,
Lemma~1.7.1 in \cite{Ross}):

\begin{lem}[Markov Inequality]\label{lem:markov}
  Let $X:\Omega\to \mathbb R$ be a nonnegative random variable on a
  sample probability space $\Omega$ with the expected value $E(X)>0$.
  Then for any $\delta>0$ we have
\[
P\big(X\ge \delta\big)\le \frac{E(X)}{\delta}.
\]
\end{lem}

\begin{lem}\label{lem:probab}

  Let $s>1$ be a fixed integer and let $A_s$ be an $s$-letter
  alphabet. Let $\Omega\subset A_s^{\ast}$ be a nonempty subset
  equipped with a discrete non-vanishing probability measure $\Pi$, so
  that $\sum_{x\in \Omega} P(\{x\})=1$. Denote $\mu(x):=P(\{x\})$ for
  any $x\in \Omega$.

  Then for any $\delta>0$ we have

\[
P\big(K_s(x)\ge -\log_2 \mu(x) -\log_2 \delta \big)=P\big( 2^{-K_s(x)}
\le \delta \mu(x)\big) \ge 1-\frac{1}{\delta}.
\]
\end{lem}

\begin{proof}
  Consider the function $X:\Omega\to \mathbb R$ defined by
  $X(x)=\frac{2^{-K_s(x)}}{\mu(x)}$.

  The $\Pi$-expected value of $X$ is

\[
E(X)=\sum_{x\in \Omega} \mu(x) \frac{2^{-K_s(x)}}{\mu(x)} \le
\sum_{x\in \Omega}2^{-K_s(x)}\le \sum_{x\in A_s^{\ast}} 2^{-K_s(x)}\le
1,
\]
where the last inequality holds by Proposition~\ref{prop:use}.

Therefore by Markov's inequality

\[
P\big(\frac{2^{-K_s(x)}}{\mu(x)} \ge \delta\big) \le
\frac{E(X)}{\delta}\le \frac{1}{\delta},
\]
and so
\[
P\big(\frac{2^{-K_s(x)}}{\mu(x)} \le \delta\big)\ge
P\big(\frac{2^{-K_s(x)}}{\mu(x)} < \delta\big)\ge 1-\frac{1}{\delta},
\]
as required.
\end{proof}

\section{Kolmogorov complexity and freely reduced words}

\begin{conv}\label{conv:1}
  Let $k>1$ and let $F=F(a_1,\dots, a_k)$. Put \[ A_{2k}:=\{a_1,\dots,
  a_k, a_1^{-1}, \dots, a_k^{-1}\}.\] As usual we identify $F$ with
  the set of all freely reduced words in $A_{2k}^{\ast}$. Thus if
  $g\in F$ then $|g|$ is the length of the unique freely reduced word
  representing $g$. For a subset $S\subseteq F$ denote by
  $\gamma(n,S)$ the number of all $x\in S$ with $|x|=n$. Similarly,
  denote by $\rho(n,S)$ the number of all $x\in S$ such that $|x|\le
  n$. Note that $\gamma(n,F)=2k(2k-1)^{n-1}$ for $n\ge 1$. Denote by
  $CR$ the set of all cyclically reduced words in $A_{2k}^{\ast}$.
  Thus $CR\subseteq F$. These notations will be fixed for the
  remainder of the paper, unless specified otherwise.
\end{conv}

It is easy to see that:
\begin{lem}\label{lem:count}\cite{KSS}
  For any $n\ge 1$ we have
\[
(2k-1)^n\le \gamma(n,CR) \le 2k (2k-1)^n.
\]
\end{lem}

Moreover, in Proposition~\ref{prop:rivin} below we will see an
explicit formula for $\gamma(n,CR)$, which we do not need for the
moment.

\begin{prop}\label{prop:main}
  Let $c\ge 1$. Denote by $Z$ the set of all cyclically reduced words
  $x$ such that
\[
C_{2k}(x)\ge -\frac{c}{2}+|x|\frac{\log_2(2k-1)}{2}.
\]

Then there is $n_0>1$ such that for any $n\ge n_0$ we have
\[
\frac{\gamma(n,Z)}{\gamma(n,CR)}\ge 1-\frac{1}{2^c}.
\]

\end{prop}

\begin{proof}
  Let $n>0$ be an integer and let $\mathcal{W}_n$ be the set of all
  cyclically reduced words of length $n$ with the uniform discrete
  probability measure $P$.  As in Lemma~\ref{lem:probab} denote
  $\mu(x):=P(\{x\})$ for any $x\in\mathcal{W}_n$.  Then by
  Lemma~\ref{lem:count} for any $x\in \mathcal{W}$ we have

\[
\frac{1}{2k} (2k-1)^{-n} \le P(\{x\})=\mu(x)=\frac{1}{\gamma(n,CR)}\le
(2k-1)^{-n}.
\]

We apply Lemma~\ref{lem:probab} with $\delta=2^c$. Hence

\begin{align*}
  &1-\frac{1}{2^c} \le P\big(2^{-K_{2k}(x)} \le 2^{c}\mu(x)\big)\le
  P\big(2^{-K_{2k}(x)}
  \le 2^{c}(2k-1)^{-n}\big)= \\
  &=P\big(-K_{2k}(x) \le c-n\log_2(2k-1))=P(K_{2k}(x)\ge
  -c+n\log_2(2k-1)\big)
\end{align*}

Recall that by Proposition~\ref{prop:use}
\[
K_{2k}(x)\le C_{2k}(x)+\log_2 C_{2k}(x)+c_0
\]
where $c_0$ is some fixed constant.  There is $n_0>1$ such that for
any word $x\in A_{2k}^{\ast}$ of length $n\ge n_0$ we have
\[
K_{2k}(x)\le 2C_{2k}(x).
\]

Therefore if $n\ge n_0$ then
\begin{align*}
  1-\frac{1}{2^c}&\le \\
  &\le P\big(K_{2k}(x)\ge -c+n\log_2(2k-1)\big)\le \\
  &\le P\big(C_{2k}(x)+\log_2 C_{2k}(x)+c_0\ge -c+n\log_2(2k-1)\big)\le\\
  &\le P\big(2C_{2k}(x)\ge -c+n\log_2(2k-1)\big),
\end{align*}
as required.
\end{proof}

\section{Genericity in free groups}\label{sect:gen}

If $b_n, b\in \mathbb R$ and $\lim_{n\to\infty} b_n=b$, we say that
the convergence is \emph{exponentially fast} if there exist $C>0$ and
$\sigma$ with $ 0<\sigma<1$ such that for all $n$ we have
\[
|b_n-b|\le C \sigma^n.
\]

\begin{defn}
  Let $S\subseteq Q\subseteq F$.

  We say that $S$ is \emph{$Q$-generic} if

\[
\lim_{n\to\infty} \frac{\rho(n,S)}{\rho(n,Q)}=1.
\]
If in addition the convergence in the above limit is exponentially
fast, we say that $S$ is \emph{exponentially $Q$-generic}.

Similarly, $S$ is called \emph{(exponentially) $Q$-negligible} if
$Q-S$ is (exponentially) $Q$-generic.
\end{defn}
Note that the union of two (exponentially) $Q$-negligible sets is
(exponentially) $Q$-negligible and the intersection of two
(exponentially) $Q$-generic sets is (exponentially) $Q$-generic.

\begin{prop}\label{prop:crit}\cite{KSS}
  The following hold:

\begin{enumerate}
\item A subset $S\subseteq F$ is exponentially $F$-negligible if and
  only if \[\lim_{n\to\infty} \frac{\gamma(n,S)}{(2k-1)^n}=0\] with
  exponentially fast convergence.
\item A subset $S\subseteq CR$ is exponentially $CR$-negligible if and
  only if \[\lim_{n\to\infty} \frac{\gamma(n,S)}{(2k-1)^n}=0\] with
  exponentially fast convergence.
\item A subset $Q\subseteq CR$ is exponentially $CR$-generic if and
  only if \[\lim_{n\to\infty} \frac{\gamma(n,Q)}{\gamma(n,CR)}=1\]
  with exponentially fast convergence.
\end{enumerate}

\end{prop}

\begin{defn}
An automorphism $\tau: F\to F$ is called a \emph{relabeling}
automorphism if the restriction $\tau|_{A_{2k}}$ is a permutation of
$A_{2k}$.
\end{defn}

\begin{conv}
For the remainder of the paper we adopt the following
convention. If $r\ge 0$ is a real number, by saying that $w$ is a word
of length $r$ we will mean that $w$ is a word of length $\lfloor r \rfloor$.
\end{conv}

\begin{lem}\label{lem:1}
  Let $0<\lambda<1/3$.  Let $\tau$ be a nontrivial relabeling
  automorphism of $F$.

  Define $S(\lambda,\tau)$ as the set of all cyclically reduced words
  $x$ such that $x$ and some cyclic permutation of $\tau(x)$ have a
  common initial segment of length $\ge \lambda |x|$.

  Then $S(\lambda,\tau)$ is exponentially $CR$-negligible.
\end{lem}
\begin{proof}
  Suppose $x\in S(\lambda,\tau)$ and $|x|=n>1$. Then there exist an
  initial segment $u$ of $x$ with $|u|=\lambda n$ and a
  cyclic permutation $\nu$ taking $\tau(x)$ to $x'$ such that $u$ is
  also an initial segment of $x'$.

\noindent{\bf Case 1.} Suppose first that $\nu$ is a trivial cyclic
permutation. Then $u$ is an initial segment of $\tau(x)$ and
$u=\tau(u)$. Since $\tau$ is a relabeling automorphism, this implies
that there is some letter $a\in A_{2k}$ such that $a^{\pm 1}$ does not
occur in $u$. Then the number of possibilities for $u$ is at most
$2k(2k-3)^{\lambda n-1}$ and the number of possibilities for $v$ is at
most $2k(2k-1)^{(1-\lambda)n-1}$.  Hence the number of all such $u$ is
at most
\[
\frac{4k^2}{(2k-1)(2k-3)} (2k-3)^{\lambda n}(2k-1)^{(1-\lambda)n}
\]
which is exponentially smaller than $(2k-1)^n$.

\noindent{\bf Case 2.} Suppose now that $\nu$ is a nontrivial cyclic
permutation, so that $\nu$ has ``translation length'' $l\ne 0 (\mod
n)$, $1\le l\le n-1$.  Thus $\tau(x)=y_1y_2$, $x'=y_2y_1$ and
$|y_2|=l$ and $|y_1|=n-l$.

The idea is that there are at least $\lambda n/6 $ letters of $x$ for
which there is no choice and which are predetermined by the rest of
$x$. Hence the number of possibilities for $x$ is exponentially
smaller than $(2k-1)^n$. There are basically two cases: when the
overlap between the positions of $u$ in $x$ and in $\tau(x)$ is small
(that is both $l$ and $-l$ are large $\mod n$) and when the  overlap
is large (that is one of $l,-l$ is small $\mod n$).

\noindent{\bf Subcase 2.A.} Assume first that $l, n-l\ge |u|/6=\lambda
n/6$, so that the overlap between the positions of $u$ in $x$ and
$\tau(x)$ has length at most $|u|/6$.

Then $y_2=uy_2'$ and $\tau(x)=y_1uy_2'$ where $|y_1|\ge |u|$.  Hence
$x=uv=uv_1\tau^{-1}(u)v_2$ where $|uv_1|=|y_1|$. We see that in this
case the segment $u'=\tau(u)$ of $x$ of length $\lambda n/6$ occurring
in the same position in $x$ as $u$ does in $\tau(x)$ is uniquely
determined (for a fixed $l$) by the rest of the word $x$. The number
of choices for $l$ is at most $n$. Given $l$ the number of choices for
$(u v_1, v_2)$ is at most $\frac{(2k)^2}{(2k-1)^2} (2k-1)^{n-\lambda
  n/6}$. Hence the number of possibilities for such $u$ is at most
\[
n\frac{(2k)^2}{(2k-1)^2} (2k-1)^{n(1-\lambda/6)}
\]
which is exponentially smaller than $(2k-1)^n$.

\noindent{\bf Subcase 2.B.} Suppose now that $0<l<|u|/6$ or $0<n-l<|u|/6$.

We will assume that $0<l<|u|$ as the other case is similar.  Thus
$x=uv$ and $\tau(x)=y_1uy_2$ with $|y_1|=l$. So the positions in which
$u$ occurs in $x$ and in $\tau(x)$ have an overlap of length $|u|-l$.
That is we can write $u=z_1u_1$ with $|u_1|=l$

Represent $|u|=m_0 l +d_0$ with $0\le d_0<l$.  Note that $m_0\ge 5$
and $d_0<l\le |u|/6=\lambda n/6$.

Now write $u$ as
\[
u=z'u_{m_0}u_{m_0-1}\dots u_1
\]
where $|u_i|=l$ for $i=1,\dots, m_0$ and $|z'|=d_0$.

Since

\begin{align*}
  x=uv=z'u_{m_0}&u_{m_0-1}\dots u_2u_1v\text{  and }\\
  \tau(x)=y_1uy_2=y_1 z' &u_{m_0}u_{m_0-1}\dots u_2u_1y_2,
\end{align*}
and $|y_1|=l$, we see that
\[
u_2=\tau(u_1), u_3=\tau(u_2)=\tau^2(u_1), \dots
u_{m_0}=\tau(u_{m_0-1})=\tau^{m_0-1}(u_1).
\]
Thus, given $l$, the words $u_1$ and $z'$ determine uniquely the rest
of the word $u$, namely the word $w=u_{m_0}\dots u_2$. Recall that
$|z'|\le l$, $|u_1|=l$ and hence
\[|w|\ge |u|-2l\ge |u|-2|u|/6=2|u|/3=2\lambda n/3.\]

Recall that $|z'|=d_0$ is determined by $l$.  So, given $l$ (for which
there are at most $n$ choices), the word $w$ is uniquely determined by
the rest of the word $x$.

Hence the number of possibilities for $x$ is at most

\[
n \frac{(2k)^2}{(2k-1)^2} (2k-1)^{n-2\lambda n/3}
\]
which is exponentially smaller than $(2k-1)^n$.

By summing up the numbers of possibilities for $x$ in the above cases
we see that

\[
\lim_{n\to\infty}\frac{\gamma(n,S(\lambda,\tau))}{(2k-1)^n} =0
\]
with exponentially fast convergence.

Hence $S(\lambda,\tau)$ is exponentially $CR$-generic by
Proposition~\ref{prop:crit}.

\end{proof}

The same type of an argument as in the proof of Lemma~\ref{lem:1}
yields:

\begin{lem}\label{lem:2}
  Let $0<\lambda<1/3$.  Let $\tau$ be a nontrivial relabeling
  automorphism of $F$.

  Define $S'(\lambda,\tau)$ as the set of all cyclically reduced words
  $x$ such that $x$ and some cyclic permutation of $\tau(x^{-1})$ have
  a common initial segment of length $\ge \lambda |x|$.  Then
  $S'(\lambda,\tau)$ is exponentially $CR$-negligible.
\end{lem}

\begin{defn}\label{defn:Y}
  Let $0<\lambda<1/3$. For a non-proper power cyclically reduced word
  $x$ let $\mathcal Y(x,\lambda)$ be the set of all $y$ satisfying one
  of the following:
\begin{enumerate}
\item the word $y$ is a cyclic permutation of $\tau(x)$ for some
  nontrivial relabeling automorphism $\tau$;
\item the word $y$ is a cyclic permutation of $\tau(x^{-1})$ for some
  (possibly trivial) trivial relabeling automorphism $\tau$;
\item the word $y$ is obtained by a nontrivial cyclic permutation of
  $x$,
\end{enumerate}
\end{defn}

\begin{lem}\label{lem:3}
  Let $0<\lambda<1/3$. Define $E(\lambda)$ as the set of all
  non-proper power cyclically reduced words $x$ such that for every
  $y\in \mathcal Y(x,\lambda)$ the lengths of the maximal common
  initial segment of $x$ and $y$ is $< \lambda |x|$.  Then
  $E(\lambda)$ is exponentially $CR$-generic.
\end{lem}
\begin{proof}
  As proved by Arzhantseva and Ol'shanskii~\cite{AO} (and easy to see
  directly by arguments similar to those used in Lemma~\ref{lem:1} and
  Lemma~\ref{lem:2}), the set of non-proper power cyclically reduced
  words $x$ whose symmetrized closures satisfy the $C'(\lambda)$ small
  cancellation condition (see \cite{LS} for definitions) is
  exponentially $CR$-generic. Since there are only finitely many
  relabeling automorphisms, the result now follows from
  Lemma~\ref{lem:1} and Lemma~\ref{lem:2} by intersecting a finite
  number of exponentially $CR$-generic sets.
\end{proof}

\begin{rem}\label{rem:dist}
  Note that by definition the set $E(\lambda)$ is closed under taking
  inverses, cyclic permutations and applying relabeling automorphisms.
  Let $M$ be the number of all (including the trivial one) relabeling
  automorphisms. Then for any $x\in E(\lambda)$ the set $Y(x,\lambda)$
  contains exactly $2M|x|-1$ distinct elements.
\end{rem}

\section{Delzant's $T$-invariant for one-relator groups}

\begin{defn}[Non-reduced $T$-invariant]
For a finite group presentation $\Pi=\langle X|R\rangle$ denote
$\displaystyle \ell_1(\Pi):=\sum_{r\in R} |r|$.

If $G$ is a finitely presentable group, define

\[
T_1(G):=\min\{ \ell_1(\Pi) : \Pi \text{ is a finite presentation
of } G\}.
\]
We call $T_1(G)$ the \emph{non-reduced $T$-invariant of
  $G$}.
\end{defn}

Obviously, for any $\Pi$ we have $\ell(\Pi)\le \ell_1(\Pi)$ and
hence for every finitely presentable group $G$ we have $T(G)\le
T_1(G)$. It turns out that under some mild assumptions there is a
similar inequality in the other direction:

\begin{lem}\label{lem:l_1}
  Let $G$ be a finitely presentable group with no elements of order two.  Then
  there exists a finite presentation $\Pi$ of $G$ such that
  $\ell(\Pi)=T(G)$ and such that every relation in $\Pi$ has length at
  least three, and therefore $T_1(G)\le \ell_1(\Pi)\le 3\ell(\Pi)=3T(G)$.

  Consequently
\[
T(G)\le T_1(G)\le 3T(G).
\]

\end{lem}
\begin{proof}
  Among all finite presentations $\Pi$ of $G$ with $\ell(\Pi)=T(G)$,
  choose a presentation $\Pi=\langle X|R\rangle$ of minimal
  $\ell_1$-length.

  We claim that every relation in $\Pi$ has length at least three.
  Clearly, the minimality assumptions on $\Pi$ imply that $\Pi$ has no
  relations of length one. Suppose $\Pi$ has a relation $r$ of length
  two. Thus $r=xy$ where $x,y\in X^{\pm 1}$.  We may assume that $y\in
  X$.

   If  $x\ne y$ in $F(X)$,  let $\Pi'$ be the presentation
  obtained from $\Pi$ by the Tietze transformation consisting of
  replacing every occurrence of $y$ in the
  relators of $R$ different from $r$ by $x^{-1}$, freely reducing the
  resulting relators  if needed, then removing the relator  $xy$ and
  removing the generator $y$ from $X$. Then $\ell(\Pi')\le
  \ell(\Pi)=T(G)$ and hence $\ell(\Pi')=T(G)$.
  By construction, $\ell_1(\Pi')<\ell_1(\Pi)$ contradicting  the
  minimality  of $\Pi$.

   If  $r=x^2$,  the
  assumption that $G$ has no elements of order two  implies that $x=1$
  in $G$. Let $\Pi''$ be the presentation obtained from $\Pi$ by
  removing the generator $x$ from $X$, removing the relation $r=x^2$
  and deleting all the occurrences of $x$ from the other relations of
  $R$ and freely reducing the results if necessary.
  We again have  $\ell(\Pi'')=\ell(\Pi)=T(G)$ and hence $\ell(\Pi'')=T(G)$.
  By construction $\ell_1(\Pi'')<\ell_1(\Pi)$,
  contradicting  the choice of $\Pi$.

  Thus every relation in $\Pi$ has length at least three, as claimed.
\end{proof}

Recall that, as specified in Convention~\ref{conv:1}, $k>1$ is a fixed
integer and $F=F(a_1,\dots, a_k)$.  As before we identify $F$ with the
set of all freely reduced word in the alphabet $A_{2k}=\{a_1,\dots,
a_k, a_1^{-1}, \dots, a_k^{-1}\}$.  For $u\in F$ we denote by $G_u$
the one-relator group $G_u:=\langle a_1,\dots, a_k | u=1\rangle$.  If
$\Pi$ is a presentation, $G(\Pi)$ denotes the group presented by
$\Pi$.

We now recall an important result about isomorphism rigidity of
generic one-relator groups that we obtained in~\cite{KSS}.

\begin{thm}\label{thm:KSS}\cite{KSS}
  Let $k>1$ be a fixed integer and $F=F(a_1,\dots, a_k)$. There exists
  an exponentially $CR$-generic set $Q_k\subseteq CR$ with the
  following properties:

\begin{enumerate}
\item There is an exponential time algorithm which, given $w\in F$,
  decides whether or not $w\in Q_k$.
\item The set $Q_k$ is closed under taking cyclic permutations,
  inverses and applying relabeling automorphisms.
\item Each $u\in Q_k$ is minimal in its $Aut(F)$-orbit, that is
  $|u|\le |\alpha(u)|$ for any $\alpha\in Aut(F)$.
\item If $r\in Q_k$ then $G_r$ is torsion-free freely indecomposable
  non-elementary word-hyperbolic group.

\item If $u\in Q_k$ and $v\in F$ are such that $|u|=|v|$ and
  $Aut(F)u=Aut(F)v$ then $v\in Q_k$ and there is a relabeling
  automorphism $\tau$ of $F$ such that $v$ is a cyclic permutation of
  $\tau(u)$.
\item Let $u\in Q_k$ and $v\in F$ be such that $|u|=|v|$. Then
  $G_u\cong G_v$ if and only if $v\in CR$ and there is a relabeling
  automorphism $\tau$ of $F$ such that $v$ is a cyclic permutation of
  $\tau(u)$ or $\tau(u)^{-1}$.
\item If $u\in Q_k$ and $v\in F$ then $G_u\cong G_v$ if and only if
  there is $\alpha\in Aut(F)$ such that $\alpha(v)=u$ or
  $\alpha(v)=u^{-1}$.
\end{enumerate}
\end{thm}

The following lemma is just the ``general enumeration argument''.

\begin{lem}\label{lem:alg}
  Let $\mathcal C$ be a a recursively enumerable class of finite
  presentations of groups. There is a partial algorithm $\Omega
  (\mathcal{C})$ which, when given a finite presentation $\Pi=\langle
  X| R\rangle$, finds a finite presentation $\Pi' \in \mathcal{C}$
  such that $G(\Pi')$ is isomorphic to $G(\Pi)$ if such a presentation
  $\Pi'$ exists.
\end{lem}

\begin{proof}
  We assume that the generating sets for $\Pi$ and for all
  presentations in $\mathcal{C}$ are initial segments of a fixed
  recursive set $\{ x_1, x_2, ...\}$ of generators.  We enumerate
  \emph{all} tuples $(\Pi_n,b,h,h')$ where
\[ \Pi' =  \langle X' | R'\rangle \in \mathcal{C}, d \in \mathbb{N}^{+},
h : X \longrightarrow F(X'), h': X': \longrightarrow F(X) \]

When such a tuple is enumerated, we then enumerate the first $d$
elements of $N' = ncl(R') \subset F(X') \text{and of} N = ncl(R)
\subset F(X)$.  We then check all of the following hold
using only the elements of $N$ and $N'$ which
have just been enumerated:

\begin{gather*}
  h(r)\in N'  \text{ for all } r\in R,\\
  h'(r)\in N \text{ for all } r\in R',\\
  h'h(x)x^{-1}\in N  \text{ for all } x\in X,\\
  hh'(x)x^{-1}\in N' \text{ for all } x\in X'.
\end{gather*}

If all of these memberships are witnessed by the elements of $N$ and
$N'$ just enumerated, then $h$ and $h'$ define mutually inverse
isomorphisms between $G(\Pi)$ and $G(\Pi')$ and we output $\Pi'$.  If
not, we go on to the next tuple.
\end{proof}

\begin{conv}
  For $0<\lambda<1/3$ denote $Q_k(\lambda)=Q_k\cap E(\lambda)$ where
  $E(\lambda)$ is as in Lemma~\ref{lem:3} and $Q_k$ is from
  Theorem~\ref{thm:KSS}. By Lemma~\ref{lem:3} and
  Theorem~\ref{thm:KSS} the set $Q_k(\lambda)$ is exponentially
  $CR$-generic.
\end{conv}

\begin{lem}\label{lem:N}
  There exists a constant $N=N(k)>0$ with the following property.  Let
  $0<\lambda<1/3$ be a rational number and let $r\in Q_k(\lambda)$ be
  a nontrivial cyclically reduced word and $G_r:=\langle a_1,\dots,
  a_k| r=1\rangle$. Thus $r$ is not a proper power and it satisfies
  the $C'(\lambda)$ small cancellation condition.

  Suppose $G_r$ can be presented by a finite presentation
\[
\Pi=\langle b_1, \dots b_m | r_1, \dots, r_t\rangle \tag{\dag}
\]
where $t\ge 1$.

Then $C_{2k}(r)\le N\ell_1(\Pi)\log_2 \ell_1(\Pi)+|r|N\lambda+N$.
\end{lem}

\begin{proof}
  We describe an algorithm $\mathcal A$, which, given a presentation
  $(\dag)$ for $G_r$ and an initial segment $u$ of $r$ of length
  $\lambda |r|$, will recover the word $r$.

  First, note that we are assuming that $(\dag)$  defines a group
  isomorphic to the $k$-generator
  one-relator group $G_r$ with  defining relator in $Q_k$. We first apply
  the algorithm $\Omega(\mathcal{C})$ from Lemma~\ref{lem:alg} with
  $\mathcal C$ the class of all $k$-generator one-relator
  presentations with defining relators from $Q_k$.  (Note that $\mathcal
  C$ is recursive by part (1) of Theorem~\ref{thm:KSS}.) This
  procedure finds some
  cyclically reduced word $v\in Q_k$ such that $(\dag)$ defines a
  group isomorphic to $G_v$.

  Thus $G_r\cong G_v$ and both $r$ and $v$ (as well as $v^{-1}$) are
  minimal cyclically reduced words from $Q_k$. By
  Theorem~\ref{thm:KSS} $|v|=|r|$ and there is a relabeling
  automorphism $\tau$ of $F$ such that $r$ is a cyclic permutation of
  $\tau(v)$ or $\tau(v)^{-1}$.

  Construct the set $B$ consisting of all words $x$ with the property
  that there is a relabeling automorphism $\tau$ of $F$ such that $x$
  is a cyclic permutation of $\tau(v)$ or $\tau(v)^{-1}$.  Thus $r\in
  B$. By Lemma~\ref{lem:3} there is a unique element of $B$ having the
  same initial segment of length $\lambda |r|$ as does $r$, namely $r$
  itself. Recall that the initial segment $u$ of $r$ of length
  $\lambda |r|$ is  part of the input for algorithm $\mathcal A$.
  Then  we list all elements of $B$ and check which one of them has
  initial segment $u$. That element is $r$.

  The algorithm $\Omega(\mathcal{C})$ is fixed.  The further input of
  $\mathcal A$, required to compute $r$, consists of the presentation
  $(\dag)$ and the initial segment $u$ of $r$ with $|u|=\lambda |r|$.
  We need to estimate the length of this input when expressed as a binary
  sequence. Put $T=\ell_1(\Pi)$. First note that in $(\dag)$ every $b_i$
  must occur in some $r_j^{\pm 1}$ since $G_r$ is a one-ended group by
  Theorem~\ref{thm:KSS} and therefore $m\le T$.

  We can now encode the presentation $(\dag)$ by writing each subscript
  $i=1,\dots, m$ for each occurrence of $b_i$ in $(\dag)$ as a binary
  integer. Using $\overline i$ to denote the binary expression for
  $i$, we replace each occurrence of $b_i$ in $(\dag)$ by
  $b \overline i$ and each occurrence of $b_i^{-1}$ by $-b\overline i$.
  Note that the bit-length of the binary
  expression $\overline i$ of $i$ is at most $\log_2 i$.  This produces an unambiguous encoding
  of $(\dag)$ as a string $W$ of length at
  most $O(T\log_2 T)$  over the six letter  alphabet
  \[
  b\quad 0\quad 1\quad -\quad , \quad |
  \] and this
  alphabet can then be block-coded into binary in the standard way.

  Since the number $k$ of generators is fixed, describing $u$ requires at most $O(|u|)$ number
  of bits.

  Hence there exist a constant $N=N(k)>0$ such that
\[
C_{2k}(r)\le NT\log_2 T+|r|N\lambda+N.
\]

\end{proof}

\begin{thm}\label{thm:main}
  Let $k>1$ be a fixed integer and $F=F(a_1,\dots, a_k)$.  For any
  $\epsilon, 0<\epsilon<1$ there is an integer $n_1>0$ and a constant
  $M=M(k,\epsilon)>0$ with the following property.

  Let $J$ be the set of all nontrivial cyclically reduced words $r$
  such that
\[
T(G_r)\log_2 T(G_r) \ge M |r|.
\]
Then for any $n\ge n_1$

\[
\frac{\gamma(n,J)}{\gamma(n,CR)}\ge 1-\epsilon.
\]
\end{thm}

\begin{proof}
  Let $N>0$ be the constant provided by Lemma~\ref{lem:N}. Choose a
  rational number $\lambda$, $0<\lambda<1/3$ so that
  $\frac{\log_2(2k-1)}{2}- N\lambda>0$.

  Let $c>0$ be an arbitrary integer. Let $n_0>1$ be the integer
  provided by Proposition~\ref{prop:main}. As in
  Proposition~\ref{prop:main} let $Z$ be the set of all cyclically
  reduced words $x$ of length $\ge n_0$ such that
\[
C_{2k}(x)\ge -\frac{c}{2}+|x|\frac{\log_2(2k-1)}{2}.
\]

Then by Proposition~\ref{prop:main} for any $n\ge n_0$ we have
\[
\frac{\gamma(n,Z)}{\gamma(n,CR)}\ge 1-\frac{1}{2^c}
\]
Since $Q_k(\lambda)$ is exponentially generic,
Proposition~\ref{prop:crit} implies that there is $n_1\ge n_0$ such
that for any $n\ge n_1$

\[
\frac{\gamma(n,Z\cap Q_k(\lambda))}{\gamma(n,CR)}\ge 1-2\frac{1}{2^c}.
\]

Now suppose $r\in Z\cap Q_k(\lambda)$ and $|r|\ge n_1$.

Then by Lemma~\ref{lem:N}
\[
-\frac{c}{2}+|r|\frac{\log_2(2k-1)}{2}\le C_{2k}(r)\le N
T_1(G_r)\log_2 T_1(G_r)+|r|N\lambda+N
\]
and hence by Lemma~\ref{lem:l_1}

\begin{gather*}
|r|(\frac{\log_2(2k-1)}{2}- N\lambda)-N-\frac{c}{2}\le N
T_1(G_r)\log_2 T_1(G_r)\le\\
3N T(G_r)\log_2 3T(G_r)=3N T(G_r)(\log_2T(G_r)+\log_2 3)\le\\
 30 N T(G_r)\log_2 T(G_r),
\end{gather*}
yielding the conclusion of the theorem.
\end{proof}

We need the following result of Rivin on the precise number of
cyclically reduced words of a given length:

\begin{prop}\label{prop:rivin}
  For any $n\ge 1$ we have
\[
\gamma(n,CR)=(2k-1)^n+1+(k-1)[1+(-1)^n].
\]
\end{prop}

Thus for a fixed $k\ge 2$ we have $\gamma(n,CR)\sim (2k-1)^n$.

The following statement is obvious:

\begin{lem}
  The number of relabeling automorphisms is $k! 2^k$.
\end{lem}

\begin{thm}
  Fix an integer $k\ge 2$.  Let $I_n$ be the number of isomorphism
  types of groups admitting a $k$-generator one-relator presentation
  where the defining relator is cyclically reduced and has length $n$.
  Then
\[
I_n\sim \frac{(2k-1)^n}{n k! 2^{k+1}}.
\]
\end{thm}

\begin{proof}
  Choose $0<\lambda<1/3$ so that $Q_k(\lambda)\subseteq CR$ is
  exponentially $CR$-generic. Recall that $Q_k(\lambda)$ is closed
  under applying inverses, cyclic permutations and relabeling
  automorphisms.

  Denote by $M=k!2^k$ the number of all relabeling automorphisms of
  $F=F(a_1,\dots, a_k)$.

  By Remark~\ref{rem:dist} for any $u\in Q_k(\lambda)$ we have $\#
  \mathcal Y(u,\lambda)=2M|u|-1$. Hence by Theorem~\ref{thm:KSS} the
  number of all $v\in Q_k(\lambda)$ with $G_u\cong G_v$ is equal to
  $2M|u|$. Therefore the set of words of length $n$ in $Q_k(\lambda)$
  defines precisely $\frac{\gamma(n, Q_k(\lambda))}{2Mn}$ isomorphism
  types of one-relator groups. Denote
  $b_n=\gamma(n,CR)-\gamma(n,Q_k(\lambda))$. Thus
  $\frac{b_n}{(2k-1)^n}\to 0$ exponentially fast as $n\to \infty$.

  Hence

\[
|I_n-\frac{\gamma(n, Q_k(\lambda))}{2Mn}|\le b_n,
\]
and so

\[
|\frac{2nMI_n}{(2k-1)^n}-\frac{\gamma(n, Q_k(\lambda))}{(2k-1)^n}|\le
\frac{2Mnb_n}{(2k-1)^n}.
\]

By $CR$-genericity of $Q_k(\lambda)$ and by Rivin's formula we have
\[
\lim_{n\to\infty} \frac{\gamma(n,
  Q_k(\lambda))}{(2k-1)^n}=\lim_{n\to\infty}\frac{\gamma(n,
  Q_k(\lambda))}{\gamma(n,CR)} \frac{\gamma(n,CR)}{(2k-1)^n}=1\cdot
1=1.
\]

Since $\lim_{n\to\infty} \frac{2Mnb_n}{(2k-1)^n}=0$, this implies
\[
\lim_{n\to\infty} \frac{2nMI_n}{(2k-1)^n}=1,
\]
and hence

\[
I_n\sim \frac{(2k-1)^n}{2Mn}=\frac{(2k-1)^n}{n k! 2^{k+1}},
\]
as required.
\end{proof}

\end{document}